\documentclass[a4paper,12pt]{article}
\usepackage[latin1]{inputenc}

\usepackage[T1]{fontenc}

\usepackage{amssymb}

\usepackage{amsmath}
\usepackage{amsthm}
\usepackage[francais,english]{babel}
\usepackage{graphicx}
\usepackage{epsfig}
\usepackage{color}
\usepackage[all]{xy}
\usepackage{latexsym}

\usepackage{geometry}

\newcommand{\graphe}{\includegraphics{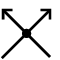}}
\newcommand{\gr}{\begin{array}{c}\graphe\end{array}}

\newcommand{\Rdeux}{\includegraphics{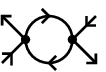}}
\newcommand{\reid}{\begin{array}{c}\Rdeux\end{array}}

\newcommand{\Rtroisd}{\includegraphics{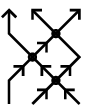}}
\newcommand{\reitd}{\begin{array}{c}\Rtroisd\end{array}}

\newcommand{\Rtroisg}{\includegraphics{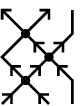}}
\newcommand{\reitg}{\begin{array}{c}\Rtroisg\end{array}}

\newcommand{\Rtroislissed}{\includegraphics{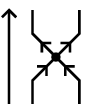}}
\newcommand{\reitlid}{\begin{array}{c}\Rtroislissed\end{array}}

\newcommand{\Rtroislisseg}{\includegraphics{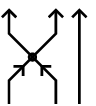}}
\newcommand{\reitlig}{\begin{array}{c}\Rtroislisseg\end{array}}

\newcommand{\lisse}{\includegraphics{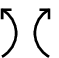}}
\newcommand{\li}{\begin{array}{c}\lisse\end{array}}

\newcommand{\gauche}{\includegraphics{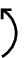}}
\newcommand{\ga}{\begin{array}{c}\gauche\end{array}}

\newcommand{\positif}{\includegraphics{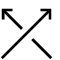}}
\newcommand{\po}{\begin{array}{c}\positif\end{array}}

\newcommand{\negatif}{\includegraphics{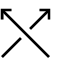}}
\newcommand{\nega}{\begin{array}{c}\negatif\end{array}}

\newcommand{\Run}{\includegraphics{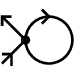}}
\newcommand{\reiu}{\begin{array}{c}\Run\end{array}}

\newcommand{\trivial}{\includegraphics{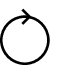}}
\newcommand{\triv}{\begin{array}{c}\trivial\end{array}}

\newcommand{\lisseinv}{\includegraphics{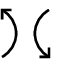}}
\newcommand{\liv}{\begin{array}{c}\lisseinv\end{array}}

\newcommand{\lisseRdeux}{\includegraphics{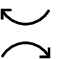}}
\newcommand{\liRde}{\begin{array}{c}\lisseRdeux\end{array}}

\newcommand{\deuxgraphe}{\includegraphics{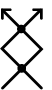}}
\newcommand{\reideux}{\begin{array}{c}\deuxgraphe\end{array}}

\newtheorem{thm}{Theorem}

\newtheorem{cor}{Corollary}

\newtheorem{lem}{Lemma}

\newtheorem{ex}{Example}

\title{Khovanov-Rozansky Graph Homology and Composition Product}

\author{Emmanuel WAGNER}

\date{ }

\begin{document}

\maketitle \abstract{In analogy with a recursive formula for the
HOMFLY-PT polynomial of links given by Jaeger, we give a recursive
formula for the graph polynomial introduced by Kauffman and Vogel.
We show how this formula extends to the Khovanov-Rozansky
graph homology.

\section*{Introduction}
The two variable HOMFLY-PT polynomial $P$ of oriented links in $\mathbb{R}^3$ is uniquely determined by its value on an unknot and by the skein relation in Figure \ref{skein}, see \cite{HOMFLY}.
\begin{figure}
\begin{center}
\scalebox{0.5}{\input{skein2.pstex_t}}
\caption{The HOMFLY-PT skein relation}
\label{skein}
\end{center}
\end{figure}
The specialization $a=q^n$ and $b=q-q^{-1}$ for a positive integer
$n$ gives a Laurent polynomial in one variable $q$. We denote this
one variable polynomial by $P_n(L)$, where $L$ is an oriented
link, or by $P_n(D)$ if $D$ is a diagram for $L$; the
normalization here is
$$P_n\left( unknot \right)=[n]_{q}=\frac{q^n-q^{-n}}{q-q^{-1}}.$$

Jaeger \cite{JAE} introduced a recursive formula for the HOMFLY-PT
polynomial. In particular, for any oriented link diagram $D$ and
any integers $m\mbox{, }n\mbox{ }\geq 1$, this formula allows a computation
of $P_{n+m}(D)$ as a sum of products $P_n(D_1)\mbox{ } P_m(D_2)$ where
$D_1$ and $D_2$ run over certain subdiagrams of $D$. Jaeger calls
this formula a {\it composition product}. 

In this paper we study finite oriented 4-valent graphs embedded in
$\mathbb{R}^2$ such that the orientation of the edges around any
vertex is as in Figure \ref{g}.
\begin{figure}
\begin{center}
\includegraphics[height=1cm]{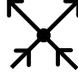}
\end{center}
\caption{Local model} \label{g}
\end{figure}
We call such graphs in $\mathbb{R}^2$ {\it regular}. We also allow components of a regular graph to be oriented circles, see Figure \ref{exemple} for an example.
\begin{figure}
\begin{center}
\includegraphics[height=1.5cm]{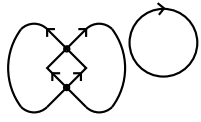}
\caption{Regular graph} \label{exemple}
\end{center}
\end{figure}
Expanding all vertices of a regular graph $\Gamma$ as in Figure
\ref{2}, Kauffman and Vogel \cite{KV} defined for any integer $n
\geq 1$ a Laurent polynomial $P_n(\Gamma)$ in one variable $q$.

\begin{figure}
\begin{eqnarray*}
P_n(\gr) & = & qP_n(\li)-q^{n}P_n(\nega)\\
         & = & q^{-1}P_n(\li)-q^{-n}P_n(\po)
\end{eqnarray*}
\caption{Graph polynomial $P_n$} \label{2}
\end{figure}

In analogy with Jaeger's composition product, we give a formula computing $P_{n+m}(\Gamma)$ as a sum of products $P_{n}(\Gamma_1)\mbox{ } P_{m}(\Gamma_2)$ where
$\Gamma_1$  and $\Gamma_2$ run over certain regular subgraphs of $\Gamma$ and $m$, $n$ $\geq 1$ are integers. More precisely, define a {\it labelling} of $\Gamma$ to be a mapping $f$ from the set of edges of $\Gamma$ to the set $\{1,2\}$ (an oriented circle is treated as an edge without vertices). We denote $\mathcal{L}(\Gamma)$ the set of labellings of $\Gamma$ that satisfy the following local condition.\\

{\bf Conservation law}: At every
vertex $v$ of $\Gamma$ the number of adjacent edges labelled by $1$ (resp.\ by $2$) directed towards $v$ is equal to the number of adjacent edges labelled by 1 (resp.\ 2) directed out of $v$.\\

Given $f\in\mathcal{L}(\Gamma)$, we can erase all edges labelled
by $2$ (resp.\ by $1$), smooth all 2-valent vertices (see Figure
\ref{9}) and obtain thus a regular graph denoted $\Gamma_{f,1}$
(resp.\ $\Gamma_{f,2}$).
\begin{figure}
\begin{center}
\includegraphics[height=0.7cm]{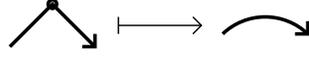}
\caption{Smoothing of a 2-valent vertex} \label{9}
\end{center}
\end{figure}

\begin{lem}For all  regular graphs $\Gamma\subset \mathbb{R}^2$ and for all integers $m$, $n$ $\geq 1,$
$$P_{n+m}(\Gamma)=\sum_{f \in \mathcal{L}(\Gamma)} q^{\sigma( \Gamma,f )}\mbox{ }P_n(\Gamma_{f,1})P_m(\Gamma_{f,2}),$$
where $\sigma( \Gamma,f )=\sigma_{m,n}( \Gamma,f )$ is an integer defined in Section 1.2.
\end{lem}

We consider the categorification of $P_n$ due to Khovanov
and Rozansky \cite{KR1}.
Given a regular graph $\Gamma \subset \mathbb{R}^2$ and an integer $n\geq 1$, Khovanov and Rozansky constructed  a $\mathbb{Z}$-graded 2-periodic chain complex $C_n(\Gamma)$ over a commutative polynomial $\mathbb{Q}$-algebra:
$$C_n(\Gamma)=\oplus_{i\in \mathbb{Z},j \in \mathbb{Z}/2\mathbb{Z}} \, C_n^{i,j}(\Gamma).$$
The differential $d$ of this complex respects the $\mathbb{Z}$-grading and increases the $\mathbb{Z}/2\mathbb{Z}$-grading by one:
$$C_n^{i,j}(\Gamma)\stackrel{d}{\rightarrow}C_n^{i,j+1}(\Gamma).$$
For all $i\in\mathbb{Z}$, $j\in\mathbb{Z}/2\mathbb{Z}$,
$$KR_n^{i,j}(\Gamma)=\mbox{Ker}(d:C_n^{i,j}(\Gamma)\rightarrow C_n^{i,j+1}(\Gamma))/\mbox{Im}(d:C_n^{i,j-1}(\Gamma)\rightarrow C_n^{i,j}(\Gamma))$$
is a finite dimensional vector space over $\mathbb{Q}$.
Set for all  $i \in \mathbb{Z}$, 
$$KR_n^i(\Gamma)=KR_n^{i,0}(\Gamma) \oplus KR_n^{i,1}(\Gamma) \mbox{ and } KR_n(\Gamma)=\oplus_{i\in \mathbb{Z},j \in \mathbb{Z}/2\mathbb{Z}}KR_n^{i,j}(\Gamma).$$ 
This construction categorifies the
graph polynomial $P_n(\Gamma)$ in the sense that
\begin{eqnarray}
P_n(\Gamma)=\sum_{i \in \mathbb{Z}} \mbox{dim}_{\mathbb{Q}} KR_n^i(\Gamma) \mbox{ }q^i.
\label{cat}
\end{eqnarray}
We will denote by  curly brackets $\{.\}$ the shift up of the $\mathbb{Z}$-grading: for $i,k\in\mathbb{Z}$ and $j\in\mathbb{Z}/2\mathbb{Z}$,  
$KR_n^{i, j}(\Gamma)\{k\}=KR_n^{i+k, j}(\Gamma)$. For $k\in \mathbb{Z}$, denote $\langle k \rangle$ the shift of the $(\mathbb{Z}/2\mathbb{Z})$-grading by $k\mbox{ (mod 2)}$.
We state now our main result.
\begin{thm}For all regular graphs $\Gamma\subset\mathbb{R}^2$, for all $m$, $n$ $\geq
1$, $i\in \mathbb{Z}$ and $j\in \mathbb{Z}/2\mathbb{Z}$,
$$KR_{n+m}^{i,j}( \Gamma)\cong \bigoplus_{\begin{array}{c}f \in \mathcal{L}(\Gamma)\\
                                                         k,l\in\mathbb{Z}, k+l+\sigma( \Gamma,f)=i\\
                                                         r,s\in\mathbb{Z}/2\mathbb{Z}, r+s=j

                                         \end{array}} KR_n^{k,r}(\Gamma_{f,1}) \otimes_{\mathbb{Q}}
KR_m^{l,s}(\Gamma_{f,2})  \{ \sigma( \Gamma,f)\}$$ where $\cong$ is
a $\mathbb{Q}$-linear isomorphism and $\sigma( \Gamma,f )=\sigma_{m,n}( \Gamma,f )$ is an integer defined in Section 1.2.

\end{thm}
This theorem yields a categorification of Lemma 1 and gives a computation of $KR_{n+m}( \Gamma)$ via
$KR_n(\Gamma_{f,1})$ and $KR_m(\Gamma_{f,2})$. We derive from Theorem 1 a direct
computational formula for $KR_{n}( \Gamma)$, see Corollary 1,
Section 3.

\-The plan of the paper is as follows. In the first section we explain all notations of Lemma 1 and Theorem 1.
The second section is devoted to the proofs of Lemma 1 and Theorem 1. In the third section, we explore consequences of Theorem 1.

\section{Preliminaries}

\subsection{Graph polynomials}

 The graph polynomial $P_n(\Gamma)\in\mathbb{Z}[q,q^{-1}]$ ($n\geq1$) of a regular graph $\Gamma\subset\mathbb{R}^2$ is defined from the relations in Figure \ref{2}.
Notice that $P_n(\Gamma)$ is preserved under ambient isotopy of
$\Gamma$ in $\mathbb{R}^2$. The polynomial $P_n$ can also be
defined as the only polynomial invariant of regular graphs invariant under ambient isotopy of graphs in $\mathbb{R}^2$ and 
satisfying the  relations in Figure \ref{3}, see \cite{KV}.
\begin{figure}
\begin{eqnarray}
P_n(\triv)=\frac{q^n-q^{-n}}{q-q^{-1}}=[n]_q \label{un}
\end{eqnarray}
\begin{eqnarray} P_n(\reiu)=[n-1]_q \mbox{ }P_n(\ga) \label{deux}
\end{eqnarray}
\begin{eqnarray}
P_n\left(\reideux \right)=[2]_q\mbox{ }P_n(\gr) \label{trois}
\end{eqnarray}
\begin{eqnarray} P_n(\reid)=P_n(\liRde)+[n-2]_q\mbox{ }P_n(\liv)
\label{quatre}
\end{eqnarray}
\begin{eqnarray}
P_n\left(\reitd\right)+P_n\left(\reitlig\right)=P_n\left(\reitg\right)+P_n\left(\reitlid\right)
\label{cinq}
\end{eqnarray}

\caption{Graph relations} \label{3}
\end{figure}
In other words, these relations are sufficient to compute $P_n(\Gamma)$ recursively.
Murakami, Ohtsuki and Yamada \cite{MOY} gave a state sum formula for $P_n(\Gamma)$ and deduced that $P_n(\Gamma)$ has only non-negative coefficients for any regular graph $\Gamma$ and any $n$ $\geq$ $1$.

\subsection{Notations}

Let $\Gamma\subset \mathbb{R}^2$ be a regular graph. We
define the {\it rotation number} of $\Gamma$. Smooth all the
vertices of $\Gamma$ as in Figure \ref{4}.
\begin{figure}
\begin{center}
\includegraphics[height=1cm]{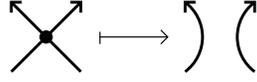}
\caption{Smoothing}
\label{4}
\end{center}
\end{figure}
This gives a disjoint union of oriented circles embedded in $\mathbb{R}^2$; we call these circles {\it Seifert circles} of $\Gamma$.
The rotation  number of $\Gamma$, denoted $r(\Gamma)$, is the sum of the signs of these circles where the sign of a  Seifert circle is $+1$ if it is oriented counterclockwise and $-1$ otherwise.

We define the interaction $\langle v|\Gamma|f \rangle\in\mathbb{Z}$ of a vertex $v$ of $\Gamma$ with a labelling $f$ as shown in Figure \ref{5}, where $1$ and $2$ are the values of $f$ on the corresponding edges.
\begin{figure}
\begin{center}
$\begin{array}{ccccc}
\begin{array}{c}\includegraphics[height=0.8cm]{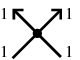}\end{array} & &
\begin{array}{c}\includegraphics[height=0.8cm]{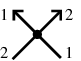}\end{array} &  &
\begin{array}{c}\includegraphics[height=0.8cm]{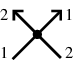}\end{array}\\
\! \langle v|\Gamma|f \rangle=0\ \! & &
\! \langle v|\Gamma|f \rangle=0\ \! & &
\! \langle v|\Gamma|f \rangle=0 \!
\end{array}$
\end{center}
\begin{center}
$\begin{array}{ccccc}
\begin{array}{c}\includegraphics[height=0.8cm]{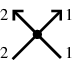}\end{array} & &
\begin{array}{c}\includegraphics[height=0.8cm]{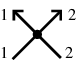}\end{array} & &
\begin{array}{c}\includegraphics[height=0.8cm]{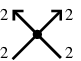}\end{array}\\
\! \langle v|\Gamma|f \rangle=1 \!& &
\! \langle v|\Gamma|f \rangle=-1 \! & &
\! \langle v|\Gamma|f \rangle=0 \!
\end{array}$
\end{center}
\caption{Definition of the interaction $\langle v|\Gamma|f \rangle$}
\label{5}
\end{figure}
Set
$\langle \Gamma|f\rangle=\sum_{v} \langle v|\Gamma|f
\rangle$, where $v$ runs over all vertices of $\Gamma$. Given integers $m$, $n$ $\geq1$, set
$$\sigma(\Gamma,f)=\sigma_{m,n}(\Gamma,f)= \langle \Gamma,f \rangle + m\, r(\Gamma_{f,1}) - n\, r(\Gamma_{f,2})\in \mathbb{Z}.$$

\section{Proofs}

\subsection{Proof of Lemma 1}
Fix $m,n\geq1$. For any regular graph $\Gamma \subset \mathbb{R}^2$, set
$$Q(\Gamma)=Q_{n+m}(\Gamma)=\sum_{f \in \mathcal{L}(\Gamma)}
\mbox{ }{q}^{\sigma(\Gamma, f)}
P_n(\Gamma_{f,1})P_m(\Gamma_{f,2})\in \mathbb{Z}[q,q^{-1}].$$In order to prove the lemma it is enough to check that $Q$  satisfies the five relations on Figure \ref{Q}.
\begin{figure}
\begin{eqnarray}
Q(\triv)=\frac{q^{n+m}-q^{-(n+m)}}{q-q^{-1}}=[n+m]_q
\label{six}
\end{eqnarray}
\begin{eqnarray}
Q(\reiu)=[n+m-1]_q\mbox{ }Q(\ga)
\label{sept}
\end{eqnarray}
\begin{eqnarray}
Q\left(\reideux \right)=[2]_q\mbox{ }Q(\gr)
\label{huit}
\end{eqnarray}
\begin{eqnarray}
Q(\reid)=Q(\liRde)+[n+m-2]_q\mbox{ }Q(\liv)
\label{neuf}
\end{eqnarray}
\begin{eqnarray}
Q\left(\reitd\right)+Q\left(\reitlig\right)=Q\left(\reitg\right)+Q\left(\reitlid\right)
\label{dix}
\end{eqnarray}
\caption{Graph relations}
\label{Q}
\end{figure}
First we verify (\ref{six}):
\begin{eqnarray*}
Q(\triv)&=& q^{-m}P_n(\triv)+q^{n}P_{m}(\triv)\\
&  =  & q^{-m}[n]_q+q^{n}[m]_q=[n+m]_q.
\end{eqnarray*}
We need to fix more notations: given a regular graph
$\Gamma$, a labelling $f_0 \in \mathcal{L}(\Gamma)$, and a subset $E_0$ of the set of edges
of $\Gamma$, set
$$Q(\Gamma_{f_0,E_0})=\sum_{f \in \mathcal{L}(\Gamma), f\vert_{E_0}={f_0}\vert_{E_0}}
\mbox{ }{q}^{\sigma(\Gamma, f)}
P_n(\Gamma_{f,1})P_m(\Gamma_{f,2}).$$
 For example,
$Q(\begin{array}{c}\includegraphics{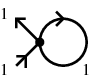}\end{array})$ is
given by a sum over all labellings whose values on the pictured
edges are 1.
We now check (\ref{sept}). %verification de la deuxième relation
We have
\begin{eqnarray}
Q(\reiu)=Q(\begin{array}{c}\includegraphics{R111.eps}\end{array})+Q(\begin{array}{c}\includegraphics{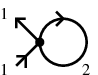}\end{array})+
Q(\begin{array}{c}\includegraphics{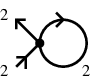}\end{array})+Q(\begin{array}{c}\includegraphics{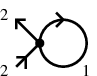}\end{array}).
\label{app}
\end{eqnarray}
Applying relation (\ref{un}) to the $P_m$-terms in $Q(\begin{array}{c}\includegraphics{R112.eps}\end{array})$, we obtain
$$Q(\begin{array}{c}\includegraphics{R112.eps}\end{array}) = q^{n-1}[m]_q\mbox{
}Q(\begin{array}{c}\includegraphics{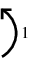}\end{array}),$$
and applying relation (\ref{deux}) to the $P_n$-terms in $Q(\begin{array}{c}\includegraphics{R111.eps}\end{array})$, we get
$$Q(\begin{array}{c}\includegraphics{R111.eps}\end{array})=
q^{-m}[n-1]_q\mbox{ }
Q(\begin{array}{c}\includegraphics{lisseg1.eps}\end{array}).$$
Similarly, we can apply relation (\ref{deux}) to the $P_m$-terms in 
$Q(\begin{array}{c}\includegraphics{R122.eps}\end{array})$ and
relation (\ref{un}) to the $P_n$-terms in
$Q(\begin{array}{c}\includegraphics{R121.eps}\end{array})$ and we
get that the right hand-side of (\ref{app}) is equal to

\begin{eqnarray*}
&   &q^{-m}[n-1]_q\mbox{ }Q(\begin{array}{c}\includegraphics{lisseg1.eps}\end{array})+q^{n-1}[m]_q\mbox{ }Q(\begin{array}{c}\includegraphics{lisseg1.eps}\end{array})\\
&   &+q^n[m-1]_q\mbox{ }Q(\begin{array}{c}\includegraphics{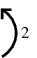}\end{array})+q^{-m+1}[n]_q\mbox{ }Q(\begin{array}{c}\includegraphics{lisseg2.eps}\end{array})\\
& = &\left(q^{-m}[n-1]_q+q^{n-1}[m]_q\right)\mbox{ }Q(\begin{array}{c}\includegraphics{lisseg1.eps}\end{array})
+\left(q^n[m-1]_q+q^{-m+1}[n]_q\right)\mbox{ }Q(\begin{array}{c}\includegraphics{lisseg2.eps}\end{array})\\
& = & [n+m-1]_q\left(Q(\begin{array}{c}\includegraphics{lisseg1.eps}\end{array})+Q(\begin{array}{c}\includegraphics{lisseg2.eps}\end{array})\right)\\
& = & [n+m-1]_q\mbox{ }Q(\ga).
\end{eqnarray*}
Hence $Q$ satisfies (\ref{sept}).\\

We check (\ref{huit}):
%verification de la troisième relation
\begin{eqnarray*}
 Q\left(\begin{array}{c}\includegraphics{R2.eps}\end{array}\right)
& = & Q\left(\begin{array}{c}\includegraphics{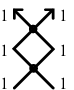}\end{array}\right)
+Q\left(\begin{array}{c}\includegraphics{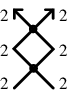}\end{array}\right)
+Q\left(\begin{array}{c}\includegraphics{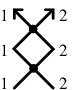}\end{array}\right)\\
&   & +Q\left(\begin{array}{c}\includegraphics{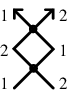}\end{array}\right)
+\mbox{ }Q\left(\begin{array}{c}\includegraphics{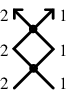}\end{array}\right)
+Q\left(\begin{array}{c}\includegraphics{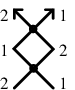}\end{array}\right)\\
&   & +Q\left(\begin{array}{c}\includegraphics{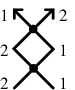}\end{array}\right)
+Q\left(\begin{array}{c}\includegraphics{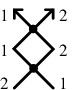}\end{array}\right)
+\mbox{ }Q\left(\begin{array}{c}\includegraphics{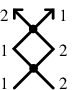}\end{array}\right)\\
&   &
+Q\left(\begin{array}{c}\includegraphics{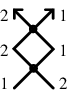}\end{array}\right)\\
%\end{eqnarray*}
%\begin{eqnarray*}
& = & Q\left(\begin{array}{c}\includegraphics{R21.eps}\end{array}\right)
+ Q\left(\begin{array}{c}\includegraphics{R22.eps}\end{array}\right)
+q^{-1} Q\left(\begin{array}{c}\includegraphics{graph12.eps}\end{array}\right)
+q Q\left(\begin{array}{c}\includegraphics{graph12.eps}\end{array}\right)\\
&  &
+\mbox{ }q Q\left(\begin{array}{c}\includegraphics{graph21.eps}\end{array}\right)
+q^{-1} Q\left(\begin{array}{c}\includegraphics{graph21.eps}\end{array}\right)
+q Q\left(\begin{array}{c}\includegraphics{graph1.eps}\end{array}\right)\\
&  & +q^{-1} Q\left(\begin{array}{c}\includegraphics{graph1.eps}\end{array}\right)
+\mbox{ }q^{-1} Q\left(\begin{array}{c}\includegraphics{graph2.eps}\end{array}\right)
+q Q\left(\begin{array}{c}\includegraphics{graph2.eps}\end{array}\right).
\end{eqnarray*}
Using the relation (\ref{trois}) for the $P_{n}$-terms and
$P_{m}$-terms,  we easily obtain (\ref{huit}).

We now check (\ref{neuf}).
%Verification de la quatrième relation
\begin{eqnarray*}
Q\left(\reid\right)
& = & Q\left(\begin{array}{c}\includegraphics{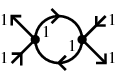}\end{array}\right)
+Q\left(\begin{array}{c}\includegraphics{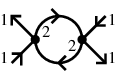}\end{array}\right)\\
&   &+Q\left(\begin{array}{c}\includegraphics{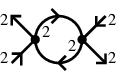}\end{array}\right)
+Q\left(\begin{array}{c}\includegraphics{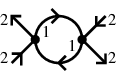}\end{array}\right)\\
&  &
+\mbox{ }Q\left(\begin{array}{c}\includegraphics{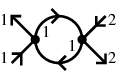}\end{array}\right)
+Q\left(\begin{array}{c}\includegraphics{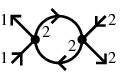}\end{array}\right)\\
&  & +Q\left(\begin{array}{c}\includegraphics{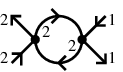}\end{array}\right)
+Q\left(\begin{array}{c}\includegraphics{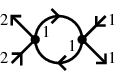}\end{array}\right)\\
&  &
+\mbox{ }Q\left(\begin{array}{c}\includegraphics{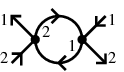}\end{array}\right)
+Q\left(\begin{array}{c}\includegraphics{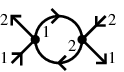}\end{array}\right).
\end{eqnarray*}
We apply  the relations (\ref{un}), (\ref{deux}), and (\ref{quatre}) to the $P_{n}$-terms and $P_{m}$-terms, so the latter expression is equal to:
\begin{eqnarray*}
&  & Q\left(\begin{array}{c}\includegraphics{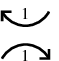}\end{array}\right)
+q^{-m}[n-2]_q\mbox{ }Q\left(\begin{array}{c}\includegraphics{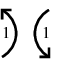}\end{array}\right)
+q^{n-2}[m]_q\mbox{ }Q\left(\begin{array}{c}\includegraphics{lisseinv1.eps}\end{array}\right)\\
&   &+Q\left(\begin{array}{c}\includegraphics{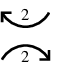}\end{array}\right)
+q^n[m-2]_q\mbox{ }Q\left(\begin{array}{c}\includegraphics{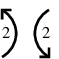}\end{array}\right)
+q^{2-m}[n]_q\mbox{ }Q\left(\begin{array}{c}\includegraphics{lisseinv2.eps}\end{array}\right)\\
&  &
+q^{1-m}[n-1]_q \mbox{ }Q\left(\begin{array}{c}\includegraphics{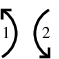}\end{array}\right)
+q^{n-1}[m-1]_q\mbox{ }Q\left(\begin{array}{c}\includegraphics{lisseinv12.eps}\end{array}\right)\\
&   &+q^{1-n}[m-1]_q\mbox{ }Q\left(\begin{array}{c}\includegraphics{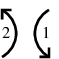}\end{array}\right)
+q^{m-1}[n-1]_q\mbox{ }Q\left(\begin{array}{c}\includegraphics{lisseinv21.eps}\end{array}\right)\\
&  &
+\mbox{ }Q\left(\begin{array}{c}\includegraphics{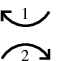}\end{array}\right)
+Q\left(\begin{array}{c}\includegraphics{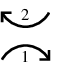}\end{array}\right)\\
& = &
Q\left(\begin{array}{c}\includegraphics{lisseR2.eps}\end{array}\right)+[n+m-2]_q\mbox{
}Q\left(\begin{array}{c}\includegraphics{lisseinv.eps}\end{array}\right).
\end{eqnarray*}

In order to prove the last relation (\ref{dix}), we put in correspondence the labellings occuring on the two
sides of (\ref{dix}) as in the following two examples:

$$Q\left(\begin{array}{c}\includegraphics{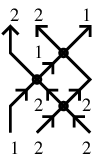}\end{array}\right)=Q\left(\begin{array}{c}\includegraphics{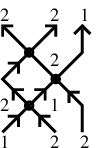}\end{array}\right)\mbox{ , }Q\left(\begin{array}{c}\includegraphics{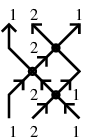}\end{array}\right)=Q\left(\begin{array}{c}\includegraphics{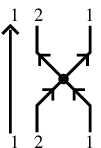}\end{array}\right).$$
There are 28 different labellings possible for $\begin{array}{c}\includegraphics{R3d.eps}\end{array}$. Among these labellings, 14 are identified with labellings of $\begin{array}{c}\includegraphics{R3g.eps}\end{array}$ as in the first example and 8 with labellings of $\begin{array}{c}\includegraphics{R3lissed.eps}\end{array}$ as in the second example. The remaining 6 labellings are involved in the following equalities (three other equalities are obtained by exchanging $1$ and $2$):

$$Q\left(\begin{array}{c}\includegraphics{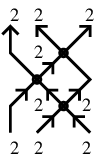}\end{array}\right)+Q\left(\begin{array}{c}\includegraphics{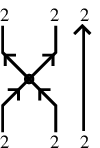}\end{array}\right)=Q\left(\begin{array}{c}\includegraphics{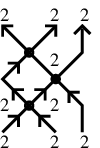}\end{array}\right)+Q\left(\begin{array}{c}\includegraphics{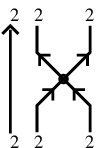}\end{array}\right),$$

$$Q\left(\begin{array}{c}\includegraphics{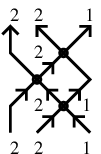}\end{array}\right)+Q\left(\begin{array}{c}\includegraphics{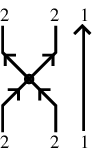}\end{array}\right)=Q\left(\begin{array}{c}\includegraphics{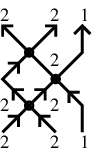}\end{array}\right),$$

$$Q\left(\begin{array}{c}\includegraphics{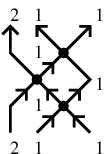}\end{array}\right)=Q\left(\begin{array}{c}\includegraphics{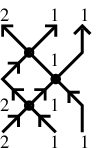}\end{array}\right)+Q\left(\begin{array}{c}\includegraphics{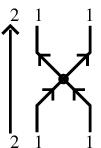}\end{array}\right).$$
In this way we obtain that $Q$ satisfies (\ref{dix}). Hence, $Q$
satisfies equations (\ref{six}-\ref{dix}) and since these
equations determine $P_{n+m}$, we conclude that
$Q=P_{n+m}$. Lemma 1 is proved.

\subsection{Proof of Theorem 1}

Since the graph polynomial has only non negative coefficients, we
have directly from Lemma 1 and Formula (\ref{cat}):
\begin{eqnarray}KR_{n+m}^{i}( \Gamma)\cong \bigoplus_{\begin{array}{c}f \in \mathcal{L}(\Gamma)\\
                                                         k,l\in\mathbb{Z}, k+l+\sigma( \Gamma,f)=i
                                         \end{array}} KR_n^{k}(\Gamma_{f,1}) \otimes_{\mathbb{Q}}
KR_m^{l}(\Gamma_{f,2})  \{ \sigma( \Gamma,f)\}.
\label{H}\end{eqnarray} 
Khovanov and Rozansky \cite{KR1}, p. 78 proved the following: for any regular graph $\Gamma$, $n \geq 1$, and $i\in\mathbb{Z}$,
$$KR_n^{i,j}(\Gamma)=0\mbox{ if } j = r(\Gamma)+1 \mbox{ (mod 2)},$$
where $r(\Gamma)$ is the rotation number defined in Section 1.2.
Furthermore, given a regular graph
$\Gamma$ and $f\in\mathcal{L}(\Gamma)$, we have
$$r(\Gamma)=r(\Gamma_{f,1})+r(\Gamma_{f,2}).$$
 Theorem 1 follows from (\ref{H}) and the latter formula.

\section{Consequences}

Given a regular graph $\Gamma$, we denote $ L (\Gamma)$ the subset of labellings of $\Gamma$ such that
$\Gamma_{f,2}$ is a disjoint union of circles. We define also
$$\mathcal{S}(\Gamma)=\{ \Gamma_{f,1}| f \in L(\Gamma)\}.$$
Given a regular subgraph $\Delta \in \mathcal{S}(\Gamma)$ of $\Gamma$, there is a unique labelling $f_{\Delta} \in L(\Gamma)$ such that $\Delta=\Gamma_{f,1}$.
For $\Delta\in \mathcal{S}(\Gamma)$, set $\beta(\Gamma,
\Delta)=\langle \Gamma, f_{\Delta} \rangle$ where $\langle \Gamma, f_{\Delta} \rangle$ is defined in Section 1.2. 
We need to fix more notations. For all $k,l\in\mathbb{Z}$,
$$\mathbb{Q}\{  k \} \langle l \rangle=\oplus_{i\in\mathbb{Z}, j\in \mathbb{Z}/2\mathbb{Z}}\mbox{ }\mathbb{Q}\{ k\} \langle l \rangle^{i,j},$$
where
$$\mathbb{Q}\{  k \} \langle l \rangle^{i,j}=\left\{\begin{array}{cc}\mathbb{Q} & \mbox{if } i=k \mbox{ and } j=l \mbox{(mod 2)},\\
                                          0 & \mbox{otherwise.}\end{array}\right.$$
We state a corollary of Theorem 1.

\begin{cor}For all regular graphs $\Gamma\subset \mathbb{R}^2$ and all integers $n\geq 2$,
$$KR_{n}(\Gamma) \cong \bigoplus_{\Delta_1 \in \mathcal{S}(\Gamma), \Delta_2 \in \mathcal{S}(\Delta_1),\dots, \Delta_{n-1} \in \mathcal{S}(\Delta_{n-2})}\mathbb{Q}\{\delta(\Delta_1,\dots,\Delta_{n-1})\} \langle r(\Gamma) \rangle$$
where
\begin{eqnarray*}
\delta(\Delta_1,\dots,\Delta_{n-1})& = &\sum_{i=0}^{n-2}\left(\beta(\Delta_i, \Delta_{i+1})+(n-i)\mbox{ }r(\Delta_{i+1})-(n-1-i)\mbox{ }r(\Delta_{i})\right)\\
                                          & = &\sum_{i=0}^{n-2}\left(\beta(\Delta_i, \Delta_{i+1})+2\mbox{ }r(\Delta_{i+1})\right) -(n-1)\mbox{ }r(\Gamma)
\end{eqnarray*}
with the convention $\Delta_0=\Gamma$.

\end{cor}
\begin{proof}
From \cite{KR1} we have
$$KR_1(\Gamma)=\left\{\begin{array}{cc}\mathbb{Q} & \mbox{if }\Gamma \mbox{ is a union of circles,}\\
                                          0 & \mbox{otherwise.}\end{array}\right.$$
The rest of the argument is just a straightforward recurrence using the case $m=1$ of Theorem 1.
\end{proof}
We illustrate Theorem 1 and Corollary 1, for $n=2$ and $m=1$ in the following example. 
\begin{ex}
\begin{eqnarray*}
KR_3\left(\begin{array}{c} \includegraphics[height=0.8cm]{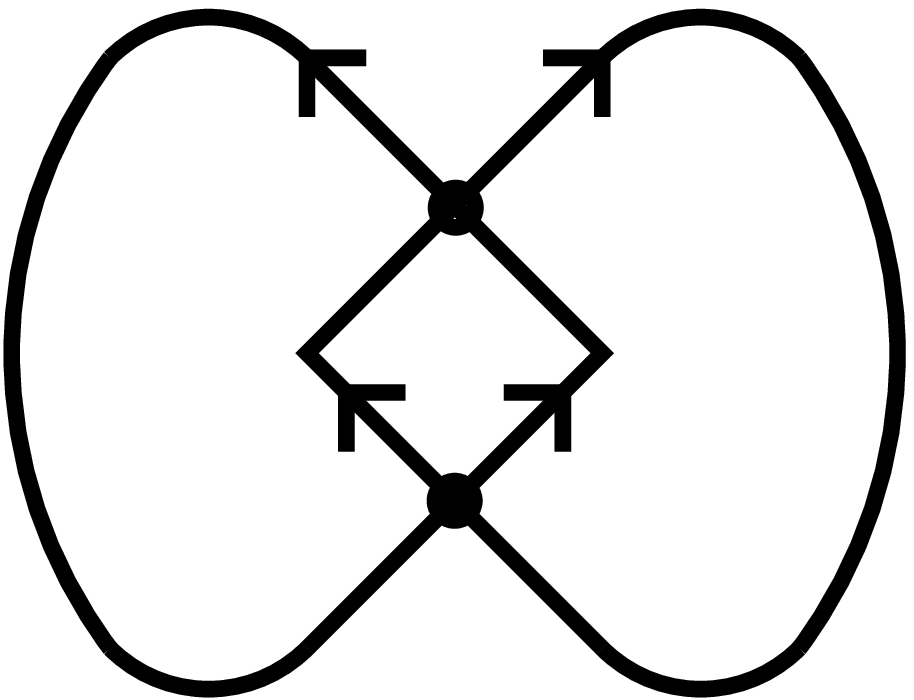}\end{array}\right) & \cong &  KR_2\left(\begin{array}{c} \includegraphics[height=0.8cm]{exe.eps}\end{array}\right) \oplus
                   KR_2(\triv)\{-3\}\langle 1 \rangle \\
&       & \oplus  KR_2(\triv)\{3\}\langle 1 \rangle \oplus  KR_2(\triv)\{1\}\langle 1 \rangle \\
&       & \oplus KR_2(\triv)\{-1\}\langle 1 \rangle \\
& \cong & \mathbb{Q}\{-2\} \oplus \mathbb{Q}\{2\} \oplus \mathbb{Q} \oplus \mathbb{Q} \oplus \mathbb{Q}\{-2\} \oplus \mathbb{Q}\{-4\}\\
&       & \oplus \mathbb{Q}\{4\} \oplus \mathbb{Q}\{2\} \oplus \mathbb{Q}\{2\} \oplus \mathbb{Q} \oplus \mathbb{Q} \oplus \mathbb{Q}\{-2\}
\end{eqnarray*}
\end{ex}
\-{\bf Remark.}
Corollary 1 suggest an equivalent but direct definition of $KR_n(\Gamma)$.
$$KR_{n}(\Gamma) = \bigoplus_{\Delta_1 \in \mathcal{S}(\Gamma), \Delta_2 \in \mathcal{S}(\Delta_1),\dots, \Delta_n \in \mathcal{S}(\Delta_{n-1})}\mathbb{Q}\{\delta(\Delta_1,\dots,\Delta_n)\} \langle r(\Gamma) \rangle.$$
%In particular, we can see that $KR_{n-1}(\Gamma)\{r(\Gamma)\}$ injects into $KR_{n}(\Gamma)$ and more generally $KR_{n-k}(\Gamma)\{kr(\Gamma)\}$ injects into $KR_{n}(\Gamma)$ for $k=1,\dots,n-1$.
Using this expression as a starting definition of $KR_n(\Gamma)$, it would be interesting to exhibit explicitly the isomorphism of Theorem 1.\\
\-{\bf Acknowledgments} This paper is done as part of the author's PhD research at the Université Louis Pasteur (Strasbourg). The author wishes to thank Prof. Turaev for his endless support.

Institut de Recherche Mathématique Avancée, C.N.R.S.-Université
Louis Pasteur, 7 rue René Descartes 67084 Strasbourg Cedex,
France.

E-mail address: wagner@math.u-strasbg.fr
\end{document}